\documentclass[12pt,diff_eng]{article}
\usepackage[english]{babel}
\usepackage{amsmath,amssymb,amscd,amsmath}
\usepackage{epsfig}
\usepackage{diff_eng}
\usepackage{graphics}
\usepackage{hyperref}
\newcommand{\tr}{^{\scriptsize\mbox{T}}}
\newcommand{\nl}{\medskip\\}
\newcommand{\disp}{\displaystyle}

\makeatletter
\def\maketitle{
	\begin{center}
		\textbf{\@title}\\
		\@author
	\end{center}
	\par
	\par}
\makeatother

\title{A Numerical-Analytical Method for Constructing Periodic Solutions of the Lorenz System}

\author{Alexander N. Pchelintsev \\
	Tambov State Technical University,\\ ul. Sovetskaya 106, Tambov, 392000, Russia \\
	e-mail: \href{mailto:pchelintsev.an@yandex.ru}{pchelintsev.an@yandex.ru}
}

\begin{document}
	\thispagestyle{empty}

\parbox{0.275\textwidth }{
\begin{picture}(130,120)
\put(67,67){{\Huge $\frac{dx}{dt}$}}
\put(40,0){\vector(0,1){70}}
\put(40,70){\line(0,1){10}}
\qbezier(40,80)(40,110)(80,110)
\qbezier(80,110)(120,110)(120,80)
\put(120,80){\line(0,-1){10}}
\qbezier(120,70)(120,35)(75,35)
\put(75,35){\vector(-1,0){5}}
\put(0,35){\line(1,0){70}}
\put(0,40){\vector(1,0){10}}
\qbezier(10,40)(35,40)(35,70)
\put(35,70){\line(0,1){10}}
\qbezier(35,80)(35,115)(80,115)
\qbezier(80,115)(125,115)(125,80)
\put(125,80){\line(0,-1){10}}
\qbezier(125,70)(125,30)(75,30)
\put(75,30){\line(-1,0){15}}
\qbezier(60,30)(45,30)(45,15)
\put(45,15){\vector(0,-1){15}}
\end{picture}
}
\hfill
\noindent \parbox{0.45\textwidth }{\footnotesize \it
\begin{center}
DIFFERENTIAL EQUATIONS\\ 
AND\\
CONTROL PROCESSES

\noindent N. 4, 2020\\

\noindent Electronic Journal, \\
reg. N ${\Phi}$C77-39410 at 15.04.2010

\noindent ISSN 1817-2172
\medskip

\noindent http://diffjournal.spbu.ru/\\

\noindent e-mail: jodiff@mail.ru

\end{center}}
\vspace*{10mm}
\large

\rightline{\underline{\it Numerical methods}}

\bigskip

	\maketitle

\noindent\textbf{Abstract.}
This article describes a method for constructing approximations to periodic solutions of dynamic
Lorenz system with classical values of the system parameters. The author obtained a system of
nonlinear algebraic equations in general form concerning of the cyclic frequency, constant terms
and amplitudes of harmonics that make up harmonic approximations to the desired solutions. The
initial approximation for the Newton method is selected, which converges to a solution describing
a periodic solution different from the equilibrium position. The results of a computational
experiment are presented. The results are verified using high-precision calculations.


\section{Introduction}
Let us consider the nonlinear system of differential equations introduced by
E. Lorenz in \cite{1}
\begin{equation}\label{eq1}
\left\{\begin{array}{l}
\dot{x}_1=\sigma(x_2-x_1),\nl
\dot{x}_2=rx_1-x_2-x_1x_3,\nl
\dot{x}_3=x_1x_2-bx_3,
\end{array}\right.
\end{equation}
where $\sigma=10$, $r=28$ and $b=8/3$ are the classical values of the system
parameters.

Let us denote by $X(t)=\left[x_1(t)\:\:x_2(t)\:\:x_3(t)\right]\tr$.
It is proved in the article \cite{1} that there exists a number $C>0$ such that
for any solution $X(t)$ of the system \eqref{eq1}, starting at time moment,
$|X(t)|<C$, and the divergence of the vector velocity field of the system \eqref{eq1}
is negative everywhere in $\mathbb R^3 $ for classical values of the system parameters.
Then \cite{1} there exists a limit set, called the Lorenz attractor, to which all
trajectories of the dynamical system are attracted when time tends to infinity. Thus
the attractor determines the behavior of the solutions of a dynamical system over
large segments of time.

W. Tucker in his work \cite{2} proved that the attractor is hyperbolic in the system
\eqref{eq1}, that is, the attractor consists of cycles everywhere dense on it along which
the near trajectories diverge exponentially. This creates their chaotic behavior.

As know \cite{3,4}, the symbolic dynamics is used to track cycles in the Lorenz system.
The region in the phase space containing the attractor is divided into a finite number of
subdomains. Denoting each partition element by a symbol, the trajectories on the attractor
passing through the corresponding regions are encoded by sequences of such symbols.
If the sequence has regularity (repeatability of groups of characters), then we have a
cycle. However, the return of trajectories in a neighborhood of its part does not mean its
closure. A critique of the results of such computational experiments can be found, for
example, in \cite{5}.

In 2004, D. Viswanath published the paper \cite{6}, in which he presented the initial
conditions and periods for three cycles in the Lorenz attractor with a high accuracy.
The calculation algorithm is based on the Lindstedt-Poincar\'e (LP) method, which
(unlike numerical integration methods) is not affected by the stability of the cycle
to which approximations are constructed.

An analysis of the Viswanath's articles \cite{6,7} showed that the author gives a general
description of the algorithm without reference to the computer implementation (in MATLAB as
indicated in his works). Moreover, it is not clear how the obtained inhomogeneous linear
system of differential equations with periodic coefficients is symbolically solved by the
LP-method. For example, this can be done for the Van der Pol equation without any special
problems.

In the article \cite{6} Viswanath showed data that can be verified by solving the
Cauchy problem with high-precision numerical methods (for example, \cite{8}), but
the details of the algorithm are not disclosed.

Therefore, it is important here to obtain the values of the initial conditions and
the period with a given accuracy, having described in detail the implementation of
the cycles search algorithm in the system \eqref{eq1}.

The goal of this article is to develop a numerical-analytical method for constructing
approximations to periodic solutions of the Lorenz system, which is simpler to
implement than the LP-method. In this case, a system of nonlinear algebraic equations
concerning of the cyclic frequency, constant terms, and amplitudes of harmonics making
up the desired solution will be obtained in general form.

\section{A Numerical-Analytical Method}

Attempts to construct approximate periodic solutions in the system \eqref{eq1} with
were made before Viswanath (for example, \cite{9}) by the method of harmonic balance,
but with low accuracy in representing real numbers, while in the article \cite{9}
initial conditions and periods of found cycles are not indicated (only drawings with cycles
are given). Now this method is actively developing in the works of \cite{10,11,12} A. Luo to
find periodic solutions of nonlinear systems of differential equations.

Next, we describe a numerical-analytical method for constructing approximations to periodic
solutions of the system \eqref{eq1}. We make for this an approximation of the phase coordinates
on the period $T$ by trigonometric polynomials in general form with an unknown cyclic frequency
$\omega$ (since we do not know the value of $T$; in the general case, it can be an irrational
number):
$$
\begin{array}{l}
\disp x_1(t)\approx\tilde{x}_1(t)=x_{1,0}+\sum_{i=1}^h\left(c_{1,i}\cos(i\omega t)+
s_{1,i}\sin(i\omega t)\right),\nl
\disp x_2(t)\approx\tilde{x}_2(t)=x_{2,0}+\sum_{i=1}^h\left(c_{2,i}\cos(i\omega t)+
s_{2,i}\sin(i\omega t)\right),\nl
\disp x_3(t)\approx\tilde{x}_3(t)=x_{3,0}+\sum_{i=1}^h\left(c_{3,i}\cos(i\omega t)+
s_{3,i}\sin(i\omega t)\right),
\end{array}
$$
where $h$ is given number of harmonics. If $i>h$, then we assume
\begin{equation}\label{eq2}
c_{1,i}=s_{1,i}=c_{2,i}=s_{2,i}=c_{3,i}=s_{3,i}=0.
\end{equation}

By the right-hand side of the system \eqref{eq1}, we compose the residuals
$$
\begin{array}{l}
\delta_1(t)=\tilde{x}'_1(t)-\sigma[\tilde{x}_2(t)-\tilde{x}_1(t)],\\
\delta_2(t)=\tilde{x}'_2(t)-[r\tilde{x}_1(t)-\tilde{x}_2(t)-\tilde{x}_1(t)\tilde{x}_3(t)],\\
\delta_3(t)=\tilde{x}'_3(t)-[\tilde{x}_1(t)\tilde{x}_2(t)-b\tilde{x}_3(t)],
\end{array}
$$
where the prime denotes the time derivative of the function. If we make calculations in an
analytical form, then for each residual you need the following:
\begin{enumerate}
	\item Differentiate by time the corresponding trigonometric polynomial;
	
	\item Where there are products of phase coordinates, multiply the corresponding
	trigonometric polynomials, converting the products of trigonometric functions
	into sums;
	
	\item Give similar terms for each function $\cos()$ and $\sin()$ with the
	corresponding argument;
	
	\item By virtue of the equalities \eqref{eq2}, to cut off the higher-order harmonics
	from the resulting residual;
	
	\item Set the resulting residual to zero, i.e., coefficients at its harmonics.
\end{enumerate}

If we put together the found algebraic equations for each residual, we obtain a still
unclosed system of nonlinear equations concerning of unknown amplitudes
$c_{1,i}$, $s_{1,i}$, $c_{2,i}$, $s_{2,i}$, $c_{3,i}$ and $s_{3,i}$ ($i=\overline{1,h}$),
constant terms $x_{1,0}$, $x_{2,0}$ and $x_{3,0}$ and the cyclic frequency $\omega$.
The number of unknown variables in the system is $3(1+2h)+1=6h+4$, but
the equations are one less.

An additional equation can be taken from the following considerations. It is known
(see \cite{4,6}) that the desired cycles intersect the plane passing through the
equilibrium positions of the system \eqref{eq1}
\begin{equation}\label{eqpos}
O_1\left(-\sqrt{b(r-1)},\,-\sqrt{b(r-1)},\,r-1\right),\:\:O_2\left(\sqrt{b(r-1)},\,
\sqrt{b(r-1)},\,r-1\right)
\end{equation}
and parallel to the plane $x_1Ox_2$ (a Poincare section).
Then the third coordinate in the initial condition for the desired cycles is equal to
$r-1$, whence $\tilde{x}_3(0)=r-1$.

Therefore the additional equation of the system has the form:
$$
x_{3,0}+\sum_{i=1}^h c_{3,i}-27=0.
$$

The author did not find in literature of other additional information on the periodic solutions
in the Lorenz system. Note that for the three cycles found by Viswanath, in the initial
condition for the third coordinate, the number 27 was taken.

Next, we give an example of a system of equations for $h=2$:
$$
\left\{
\begin{aligned}
\omega s_{1,1}-10c_{2,1}+10c_{1,1}&=0,\\
-10s_{2,1}+10s_{1,1}-c_{1,1}\omega&=0,\\
2\omega s_{1,2}-10c_{2,2}+10c_{1,2}&=0,\\
-10s_{2,2}+10s_{1,2}-2c_{1,2}\omega&=0,\\
10x_{1,0}-10x_{2,0}&=0,\\
c_{1,1}x_{3,0}+c_{3,1}x_{1,0}+\dfrac{s_{1,1}s_{3,2}}{2}+\dfrac{s_{1,2}s_{3,1}}{2}+\omega s_{2,1}+\dfrac{c_{1,1}c_{3,2}}{2}+\dfrac{c_{1,2}c_{3,1}}{2}+c_{2,1}-28c_{1,1}&=0,\\
s_{1,1}x_{3,0}+s_{3,1}x_{1,0}+\dfrac{c_{1,1}s_{3,2}}{2}-\dfrac{c_{1,2}s_{3,1}}{2}+s_{2,1}+\dfrac{c_{3,1}s_{1,2}}{2}-\dfrac{c_{3,2}s_{1,1}}{2}-28s_{1,1}-c_{2,1}\omega&=0,\\
c_{1,2}x_{3,0}+c_{3,2}x_{1,0}-\dfrac{s_{1,1}s_{3,1}}{2}+2\omega s_{2,2}+\dfrac{c_{1,1}c_{3,1}}{2}+c_{2,2}-28c_{1,2}&=0,\\
s_{1,2}x_{3,0}+s_{3,2}x_{1,0}+\dfrac{c_{1,1}s_{3,1}}{2}+s_{2,2}-28s_{1,2}+\dfrac{c_{3,1}s_{1,1}}{2}-2c_{2,2}\omega&=0,\\
x_{1,0}x_{3,0}+x_{2,0}-28x_{1,0}+\dfrac{s_{1,2}s_{3,2}}{2}+\dfrac{s_{1,1}s_{3,1}}{2}+\dfrac{c_{1,2}c_{3,2}}{2}+\dfrac{c_{1,1}c_{3,1}}{2}&=0,\\
-c_{1,1}x_{2,0}-c_{2,1}x_{1,0}+\omega s_{3,1}-\dfrac{s_{1,1}s_{2,2}}{2}-\dfrac{s_{1,2}s_{2,1}}{2}+\dfrac{8c_{3,1}}{3}-\dfrac{c_{1,1}c_{2,2}}{2}-\dfrac{c_{1,2}c_{2,1}}{2}&=0,\\
-s_{1,1}x_{2,0}-s_{2,1}x_{1,0}+\dfrac{8s_{3,1}}{3}-\dfrac{c_{1,1}s_{2,2}}{2}+\dfrac{c_{1,2}s_{2,1}}{2}-\dfrac{c_{2,1}s_{1,2}}{2}+\dfrac{c_{2,2}s_{1,1}}{2}-c_{3,1}\omega&=0,\\
-c_{1,2}x_{2,0}-c_{2,2}x_{1,0}+2\omega s_{3,2}+\dfrac{s_{1,1}s_{2,1}}{2}+\dfrac{8c_{3,2}}{3}-\dfrac{c_{1,1}c_{2,1}}{2}&=0,\\
-s_{1,2}x_{2,0}-s_{2,2}x_{1,0}+\dfrac{8s_{3,2}}{3}-\dfrac{c_{1,1}s_{2,1}}{2}-\dfrac{c_{2,1}s_{1,1}}{2}-2c_{3,2}\omega&=0,\\
\dfrac{8x_{3,0}}{3}-x_{1,0}x_{2,0}-\dfrac{s_{1,2}s_{2,2}}{2}-\dfrac{s_{1,1}s_{2,1}}{2}-\dfrac{c_{1,2}c_{2,2}}{2}-\dfrac{c_{1,1}c_{2,1}}{2}&=0,\\
x_{3,0}+c_{3,1}+c_{3,2}-27&=0.
\end{aligned}
\right.
$$

Note that for any $h$ a similar system has solutions
$$
\begin{array}{c}
\displaystyle x_{1,0}=x_{2,0}=\pm\sqrt{b(r-1)},\:x_{3,0}=r-1,\:c_{k,i}=0,\:s_{k,i}=0,\\
\omega\:\mbox{is any number},\:\,k=\overline{1,3},\:i=\overline{1,h},
\end{array}
$$
corresponding to the equilibrium positions \eqref{eqpos}.

Therefore the resulting nonlinear system of algebraic equations has a non-unique solution.
To find its approximate solutions, we will use the Newton numerical method, whose a convergence
to the desired solution (i.e., describing a periodic solution of the system \eqref{eq1} different
from its the equilibrium positions) depends on the choice of the initial approximation.

\section{The Symbolic Computations to Obtain the System of Algebraic Equations}

Thus, to obtain an approximation to the periodic solution, we must obtain a nonlinear system
concerning of unknown decomposition coefficients and frequencies. As shown in the previous section,
even for two harmonics, the system has a bulky form. Therefore, we consider the algorithm
for performing symbolic calculations to obtain it.

When developing software \cite{13}, the Maxima math package (a computer algebra system) was
chosen. The program for obtaining the amplitudes and constant terms of the residuals for $h=2$ is
presented below.

\begin{verbatim}
/* [wxMaxima batch file version 1] [ DO NOT EDIT BY HAND! ]*/
/* [wxMaxima: input   start ] */
display2d:false$
x1:x10+c1c1*cos(1*omega*t)+s1c1*sin(1*omega*t)+
c1c2*cos(2*omega*t)+s1c2*sin(2*omega*t)$
x2:x20+c2c1*cos(1*omega*t)+s2c1*sin(1*omega*t)+
c2c2*cos(2*omega*t)+s2c2*sin(2*omega*t)$
x3:x30+c3c1*cos(1*omega*t)+s3c1*sin(1*omega*t)+
c3c2*cos(2*omega*t)+s3c2*sin(2*omega*t)$
assume(omega > 0)$
delta1:trigreduce(diff(x1,t)-(10*(x2-x1)),t)$
delta2:trigreduce(diff(x2,t)-(28*x1-x2-x1*x3),t)$
delta3:trigreduce(diff(x3,t)-(x1*x2-8/3*x3),t)$
expand(diff(delta1,cos(1*omega*t)));
expand(diff(delta1,sin(1*omega*t)));
expand(diff(delta1,cos(2*omega*t)));
expand(diff(delta1,sin(2*omega*t)));
expand(integrate(delta1,t,0,2*%pi/omega)*omega/(2*%pi));
expand(diff(delta2,cos(1*omega*t)));
expand(diff(delta2,sin(1*omega*t)));
expand(diff(delta2,cos(2*omega*t)));
expand(diff(delta2,sin(2*omega*t)));
expand(integrate(delta2,t,0,2*%pi/omega)*omega/(2*%pi));
expand(diff(delta3,cos(1*omega*t)));
expand(diff(delta3,sin(1*omega*t)));
expand(diff(delta3,cos(2*omega*t)));
expand(diff(delta3,sin(2*omega*t)));
expand(integrate(delta3,t,0,2*%pi/omega)*omega/(2*%pi));
/* [wxMaxima: input   end   ] */
\end{verbatim}

The expression \verb|display2d:false$| turns off multi-line drawing of fractions,
degrees, etc. The sign \verb|$| allows to calculate the result of an expression,
but not display it (instead of \verb|;|). The function \verb|trigreduce(expression,t)|
collapses all products of trigonometric functions concerning of the variable $t$ in
a combination of sums. Differentiation of residuals according to harmonic functions is
necessary to obtain the corresponding amplitudes. The function \verb|expand(expression)|
expands brackets (performs multiplication, exponentiation, leads similar terms).

To find the constant terms of the residuals, their integration over the period is
applied, i.e. the constant term of the $k$-residual is
$$
\dfrac{\displaystyle\omega\int_{0}^{\frac{2\pi}{\omega}}\delta_k(t)dt}{2\pi}.
$$

So that during symbolic integration the package does not ask a question about the sign
of the frequency, a command is given \verb|assume(omega > 0)$|.

A file with package commands is generated similarly for any number of $h$ harmonics
by a computer program written in C++ \cite{13}. After executing this program, the
package will output symbolic expressions to the console for the left side of the system
of algebraic equations, which will be solved in it by the Newton method.

Note that the most time-consuming operation here is symbolic integration.
For example, for 120 harmonics, the system formation time is more than 2 days.
We can here parallelize the computational process on three computers, but this will
not have a significant effect. Therefore, a system of algebraic equations must be formed
immediately. Next, we get a general form of this system. Note that when solving the system
of nonlinear equations by the Newton method, the Jacobi matrix for the left side of the
system does not invert. The Maxima package uses LU decomposition to solve a system of
linear equations at each iteration of the method.

\section{General Form of the System of Algebraic Equations}

Since the right-hand side of the \eqref{eq1} system contains nonlinearities in the form
of products of phase coordinates, let us obtain relations expressing the coefficients of
trigonometric polynomials obtained by multiplying the approximations
$\tilde{x}_1(t)\tilde{x}_3(t)$ and $\tilde{x}_1(t)\tilde{x}_2(t)$.

We consider two functions $f(t)$ and $F(t)$ represented by Fourier series
$$
\begin{array}{c}
\disp f(t)=a_0+\sum_{i=1}^\infty\left(a_i\cos(i\omega t)+b_i\sin(i\omega t)\right),\nl
\disp F(t)=A_0+\sum_{i=1}^\infty\left(A_i\cos(i\omega t)+B_i\sin(i\omega t)\right).
\end{array}
$$

Let
$$
f(t)F(t)=\alpha_0+\sum_{i=1}^\infty\left(\alpha_i\cos(i\omega t)+
\beta_i\sin(i\omega t)\right).
$$

Following the book \cite[pp. 123-125]{14}, we have the relations:
$$
\alpha_0=a_0 A_0+\dfrac{1}{2}\sum_{m=1}^\infty\left(a_m A_m+b_m B_m\right),
$$
\begin{equation}\label{eq3}
\alpha_i=a_0 A_i+\dfrac{1}{2}\sum_{m=1}^\infty\left(a_m(A_{m+i}+A_{m-i})+
b_m(B_{m+i}+B_{m-i})\right),
\end{equation}
\begin{equation}\label{eq4}
\beta_i=a_0 B_i+\dfrac{1}{2}\sum_{m=1}^\infty\left(a_m(B_{m+i}-B_{m-i})-
b_m(A_{m+i}-A_{m-i})\right).
\end{equation}

We assume that for $i>h$
$$
a_i=b_i=A_i=B_i=0.
$$

Since for our problem we find for an approximation up to and including the $h$-harmonic,
we zero all the amplitudes in the product for $i>h $, i.e.
$$
\alpha_i=\beta_i=0.
$$

Thus, we pass from the product of series to the product of trigonometric polynomials.
Also in the relations \eqref{eq3} and \eqref{eq4} we will assume \cite [p. 124]{14} that
$$
A_{m-i}=A_{i-m},\:\:B_{m-i}=-B_{i-m},\:\:B_0=0.
$$

Then we get
$$
\alpha_0=a_0 A_0+\dfrac{1}{2}\sum_{m=1}^h\left(a_m A_m+b_m B_m\right),
$$
\begin{align*}
\alpha_i&=a_0 A_i+\dfrac{1}{2}\sum_{m=1}^\infty a_m A_{m+i}+
\dfrac{1}{2}\sum_{m=1}^\infty a_m A_{m-i}+
\dfrac{1}{2}\sum_{m=1}^\infty b_m B_{m+i}+
\dfrac{1}{2}\sum_{m=1}^\infty b_m B_{m-i}=\\
&=a_0 A_i+\dfrac{1}{2}\sum_{m=1}^{h-i} a_m A_{m+i}+\dfrac{1}{2}a_i A_0+
\dfrac{1}{2}\sum_{m=1}^{i-1} a_m A_{i-m}+
\dfrac{1}{2}\sum_{m=i+1}^h a_m A_{m-i}+\\
&+\dfrac{1}{2}\sum_{m=1}^h b_m B_{m+i}+\dfrac{1}{2}b_i B_0-
\dfrac{1}{2}\sum_{m=1}^{i-1} b_m B_{i-m}+
\dfrac{1}{2}\sum_{m=i+1}^h b_m B_{m-i}=\\
&=a_0 A_i+a_i A_0+\dfrac{1}{2}\sum_{m=1}^{h-i}\left(a_m A_{m+i}+b_m B_{m+i}\right)+
\dfrac{1}{2}\sum_{m=1}^{i-1}\left(a_m A_{i-m}-b_m B_{i-m}\right)+\\
&+\dfrac{1}{2}\sum_{m=i+1}^{h}\left(a_m A_{m-i}+b_m B_{m-i}\right),
\end{align*}
\begin{align*}
\beta_i&=a_0 B_i+\dfrac{1}{2}\sum_{m=1}^\infty a_m B_{m+i}-
\dfrac{1}{2}\sum_{m=1}^\infty a_m B_{m-i}-
\dfrac{1}{2}\sum_{m=1}^\infty b_m A_{m+i}+
\dfrac{1}{2}\sum_{m=1}^\infty b_m A_{m-i}=\\
&=a_0 B_i+\dfrac{1}{2}\sum_{m=1}^{h-i} a_m B_{m+i}+
\dfrac{1}{2}\sum_{m=1}^{i-1} a_m B_{i-m}-
\dfrac{1}{2}\sum_{m=i+1}^h a_m B_{m-i}-\\
&-\dfrac{1}{2}\sum_{m=1}^{h-i} b_m A_{m+i}+b_i A_0+
\dfrac{1}{2}\sum_{m=1}^{i-1} b_m A_{i-m}+
\dfrac{1}{2}\sum_{m=i+1}^h b_m A_{m-i}=\\
&=a_0 B_i+b_i A_0+\dfrac{1}{2}\sum_{m=1}^{h-i}\left(a_m B_{m+i}-b_m A_{m+i}\right)+
\dfrac{1}{2}\sum_{m=1}^{i-1}\left(a_m B_{i-m}+b_m A_{i-m}\right)+\\
&+\dfrac{1}{2}\sum_{m=i+1}^{h}\left(-a_m B_{m-i}+b_m A_{m-i}\right).
\end{align*}

Applying the obtained formulas to calculate the products of trigonometric polynomials
to the residuals, we can write the equations for the $i$-th harmonics ($i=\overline{1,h}$
is the number of harmonics, $k=\overline{1,3}$ is residual number):
{\flushleft $k=1$:}
$$
\begin{array}{r}
i\omega s_{1,i}-10c_{2,i}+10c_{1,i}=0,\nl
-i\omega c_{1,i}-10s_{2,i}+10s_{1,i}=0,
\end{array}
$$
the equation corresponding to the constant term for the first residual is
$$
x_{1,0}-x_{2,0}=0,
$$
{\flushleft $k=2$:}
\begin{align*}
i\omega s_{2,i}-28c_{1,i}+c_{2,i}+x_{1,0}c_{3,i}+c_{1,i}x_{3,0}
&+\dfrac{1}{2}\sum_{m=1}^{h-i}\left(c_{1,m}c_{3,m+i}+s_{1,m}s_{3,m+i}\right)+\\
&+\dfrac{1}{2}\sum_{m=1}^{i-1}\left(c_{1,m}c_{3,i-m}-s_{1,m}s_{3,i-m}\right)+\\
&+\dfrac{1}{2}\sum_{m=i+1}^{h}\left(c_{1,m}c_{3,m-i}+s_{1,m}s_{3,m-i}\right)=0,
\end{align*}
\begin{align*}
-i\omega c_{2,i}-28s_{1,i}+s_{2,i}+x_{1,0}s_{3,i}+s_{1,i}x_{3,0}
&+\dfrac{1}{2}\sum_{m=1}^{h-i}\left(c_{1,m}s_{3,m+i}-s_{1,m}c_{3,m+i}\right)+\\
&+\dfrac{1}{2}\sum_{m=1}^{i-1}\left(c_{1,m}s_{3,i-m}+s_{1,m}c_{3,i-m}\right)+\\
&+\dfrac{1}{2}\sum_{m=i+1}^{h}\left(-c_{1,m}s_{3,m-i}+s_{1,m}c_{3,m-i}\right)=0,
\end{align*}
the equation corresponding to the constant term for the second residual is
$$
-28x_{1,0}+x_{2,0}+x_{1,0}x_{3,0}+\dfrac{1}{2}\sum_{m=1}^h\left(c_{1,m}c_{3,m}+
s_{1,m}s_{3,m}\right)=0,
$$
{\flushleft $k=3$:}
\begin{align*}
i\omega s_{3,i}-x_{1,0}c_{2,i}-c_{1,i}x_{2,0}
&-\dfrac{1}{2}\sum_{m=1}^{h-i}\left(c_{1,m}c_{2,m+i}+s_{1,m}s_{2,m+i}\right)-\\
&-\dfrac{1}{2}\sum_{m=1}^{i-1}\left(c_{1,m}c_{2,i-m}-s_{1,m}s_{2,i-m}\right)-\\
&-\dfrac{1}{2}\sum_{m=i+1}^{h}\left(c_{1,m}c_{2,m-i}+s_{1,m}s_{2,m-i}\right)+
\dfrac{8}{3}c_{3,i}=0,
\end{align*}
\begin{align*}
-i\omega c_{3,i}-x_{1,0}s_{2,i}-s_{1,i}x_{2,0}
&-\dfrac{1}{2}\sum_{m=1}^{h-i}\left(c_{1,m}s_{2,m+i}-s_{1,m}c_{2,m+i}\right)-\\
&-\dfrac{1}{2}\sum_{m=1}^{i-1}\left(c_{1,m}s_{2,i-m}+s_{1,m}c_{2,i-m}\right)-\\
&-\dfrac{1}{2}\sum_{m=i+1}^{h}\left(-c_{1,m}s_{2,m-i}+s_{1,m}c_{2,m-i}\right)+
\dfrac{8}{3}s_{3,i}=0,
\end{align*}
the equation corresponding to the constant term for the third residual is
$$
-x_{1,0}x_{2,0}-\dfrac{1}{2}\sum_{m=1}^h\left(c_{1,m}c_{2,m}+
s_{1,m}s_{2,m}\right)+\dfrac{8}{3}x_{3,0}=0,
$$
the additional system equation is
$$
x_{3,0}+\sum_{i=1}^h c_{3,i}-27=0.
$$

\section{The Results of the Computational Experiment}

As a result of numerous computational experiments, the initial approximation
was chosen for the cyclic frequency, constant terms, and amplitudes at
$h=h_1=5$:
$$
\begin{array}{c}
\omega=4,\:\:x_{1,0}=x_{2,0}=x_{3,0}=0,\:\:c_{1,i}=-1,\:i=\overline{1,5},\\
s_{1,j}=0,\:j=1,3,4,5,\:\:s_{1,2}=1.
\end{array}
$$

This result is remarkable in that the Newton method converges to a solution
different from the equilibrium positions. Therefore, to improve the accuracy of
the approximate periodic solution, we consider a system of algebraic equations
for the value of $h$ equal to some $h_2>h_1$.

The obtained numerical solution of the system for $h=h_1$ is taken as the initial
approximation for amplitudes with indices $i\le h_1$ for a system with $h=h_2$,
and the values of the initial approximation for amplitudes with indices $i>h_1$
are assumed to be zero.

\begin{table}[!]
	\caption{The amplitudes of harmonics for $\tilde{x}_1(t)$, $x_{1,0}=0$.}
	\centering\vspace{10pt}
	{\begin{tabular}{|c|c|c|}
			\hline
			$i$ & $c_{1,i}$ & $s_{1,i}$\\\hline
			1 & $-5.780478259196228$ & $8.56017654325353$\\\hline
			2 & 0 & 0 \\\hline
			3 & $3.160762628380509$ & $2.239212141102876$\\\hline
			4 & 0 & 0 \\\hline
			5 & $0.6958870387616096$ & $-0.7979388979225431$\\\hline
			6 & 0 & 0 \\\hline
			7 & $-0.1891992374027477$ & $-0.1864921358925765$\\\hline
			8 & 0 & 0 \\\hline
			9 & $-0.04770429623010056$ & $0.04554044367245914$\\\hline
			10 & 0 & 0 \\\hline
			11 & $0.01112322884679491$ & $0.01209138588669679$\\\hline
			12 & 0 & 0 \\\hline
			13 & $0.003061207095371694$ & $-0.002735092350544739$\\\hline
			14 & 0 & 0 \\\hline
			15 & $-6.744578887916229\cdot10^{-4}$ & $-7.748319471034087\cdot10^{-4}$\\\hline
			16 & 0 & 0 \\\hline
			17 & $-1.960718247379475\cdot10^{-4}$ & $1.665584161919807\cdot10^{-4}$\\\hline
			18 & 0 & 0 \\\hline
			19 & $4.116738805347028\cdot10^{-5}$ & $4.960493476144467\cdot10^{-5}$\\\hline
			20 & 0 & 0 \\\hline
			21 & $1.254757391175977\cdot10^{-5}$ & $-1.018054283421179\cdot10^{-5}$\\\hline
			22 & 0 & 0 \\\hline
			23 & $-2.518375902000733\cdot10^{-6}$ & $-3.173486439630506\cdot10^{-6}$\\\hline
			24 & 0 & 0 \\\hline
			25 & $-8.025338211960923\cdot10^{-7}$ & $6.230623750431923\cdot10^{-7}$\\\hline
			26 & 0 & 0 \\\hline
			27 & $1.541534734542893\cdot10^{-7}$ & $2.0292802821633\cdot10^{-7}$\\\hline
			28 & 0 & 0 \\\hline
			29 & $5.130649139299358\cdot10^{-8}$ & $-3.813725452268523\cdot10^{-8}$\\\hline
			30 & 0 & 0 \\\hline
			31 & $-9.43393531993558\cdot10^{-9}$ & $-1.297038481588497\cdot10^{-8}$\\\hline
			32 & 0 & 0 \\\hline
			33 & $-3.278552746800046\cdot10^{-9}$ & $2.333260259021725\cdot10^{-9}$\\\hline
			34 & 0 & 0 \\\hline
			35 & $5.76957885768651\cdot10^{-10}$ & $8.28626640138045\cdot10^{-10}$\\\hline
		\end{tabular}
	}
\end{table}

\begin{table}[!]
	\caption{The amplitudes of harmonics for $\tilde{x}_2(t)$, $x_{2,0}=0$.}
	\centering\vspace{10pt}
	{\begin{tabular}{|c|c|c|}
			\hline
			$i$ & $c_{2,i}$ & $s_{2,i}$\\\hline
			1 & $-2.32972926505593$ & $10.89038310357172$\\\hline
			2 & 0 & 0 \\\hline
			3 & $5.86875317198698$ & $-1.5832552129833$\\\hline
			4 & 0 & 0 \\\hline
			5 & $-0.9124249133801483$ & $-2.200556873678218$\\\hline
			6 & 0 & 0 \\\hline
			7 & $-0.7154457265566421$ & $0.3473932955614448$\\\hline
			8 & 0 & 0 \\\hline
			9 & $0.1175186702136983$ & $0.2186139734768588$\\\hline
			10 & 0 & 0 \\\hline
			11 & $0.06473984670858603$ & $-0.03723215039412078$\\\hline
			12 & 0 & 0 \\\hline
			13 & $-0.01127208646321726$ & $-0.01877739524860192$\\\hline
			14 & 0 & 0 \\\hline
			15 & $-0.005359671824365359$ & $0.003303445299126894$\\\hline
			16 & 0 & 0 \\\hline
			17 & $9.453499475830811\cdot10^{-4}$ & $0.001510235036151227$\\\hline
			18 & 0 & 0 \\\hline
			19 & $4.211022386354685\cdot10^{-4}$ & $-2.657049331814368\cdot10^{-4}$\\\hline
			20 & 0 & 0 \\\hline
			21 & $-7.363528144366622\cdot10^{-5}$ & $-1.164013765469982\cdot10^{-4}$\\\hline
			22 & 0 & 0 \\\hline
			23 & $-3.19419300699788\cdot10^{-5}$ & $2.017609175377016\cdot10^{-5}$\\\hline
			24 & 0 & 0 \\\hline
			25 & $5.47663534401654\cdot10^{-6}$ & $8.710929378319451\cdot10^{-6}$\\\hline
			26 & 0 & 0 \\\hline
			27 & $2.362852034076972\cdot10^{-6}$ & $-1.474901091428546\cdot10^{-6}$\\\hline
			28 & 0 & 0 \\\hline
			29 & $-3.94532524722541\cdot10^{-7}$ & $-6.379296603810031\cdot10^{-7}$\\\hline
			30 & 0 & 0 \\\hline
			31 & $-1.715198229248314\cdot10^{-7}$ & $1.049218598356554\cdot10^{-7}$\\\hline
			32 & 0 & 0 \\\hline
			33 & $2.776045093375681\cdot10^{-8}$ & $4.59473450493284\cdot10^{-8}$\\\hline
			34 & 0 & 0 \\\hline
			35 & $1.22681173575872\cdot10^{-8}$ & $-7.31171826830086\cdot10^{-9}$\\\hline
		\end{tabular}
	}
\end{table}

\begin{table}[!]
	\caption{The amplitudes of harmonics for $\tilde{x}_3(t)$, $x_{3,0}=23.04210397942006$.}
	\centering\vspace{10pt}
	{\begin{tabular}{|c|c|c|}
			\hline
			$i$ & $c_{3,i}$ & $s_{3,i}$\\\hline
			1 & 0 & 0 \\\hline
			2 & $7.568410271550653$ & $-9.50386584559212$\\\hline
			3 & 0 & 0 \\\hline
			4 & $-3.555327211552558$ & $-1.844710563805469$\\\hline
			5 & 0 & 0 \\\hline
			6 & $-0.4741220131932616$ & $1.279043179069961$\\\hline
			7 & 0 & 0 \\\hline
			8 & $0.4227292179138024$ & $0.1274574086305204$\\\hline
			9 & 0 & 0 \\\hline
			10 & $0.03498415351761577$ & $-0.1315337800809524$\\\hline
			11 & 0 & 0 \\\hline
			12 & $-0.03934013541135439$ & $-0.009645786231708874$\\\hline
			13 & 0 & 0 \\\hline
			14 & $-0.002660052258813564$ & $0.01145537653603837$\\\hline
			15 & 0 & 0 \\\hline
			16 & $0.003271688724557337$ & $7.33752523103949\cdot10^{-4}$\\\hline
			17 & 0 & 0 \\\hline
			18 & $2.024982256871223\cdot10^{-4}$ & $-9.206266886554897\cdot10^{-4}$\\\hline
			19 & 0 & 0 \\\hline
			20 & $-2.560063570343799\cdot10^{-4}$ & $-5.58964460662525\cdot10^{-5}$\\\hline
			21 & 0 & 0 \\\hline
			22 & $-1.542436654918173\cdot10^{-5}$ & $7.050327849098175\cdot10^{-5}$\\\hline
			23 & 0 & 0 \\\hline
			24 & $1.926014222030195\cdot10^{-5}$ & $4.25261452471065\cdot10^{-6}$\\\hline
			25 & 0 & 0 \\\hline
			26 & $1.170939944189529\cdot10^{-6}$ & $-5.225643926851625\cdot10^{-6}$\\\hline
			27 & 0 & 0 \\\hline
			28 & $-1.409525591131397\cdot10^{-6}$ & $-3.21879984959824\cdot10^{-7}$\\\hline
			29 & 0 & 0 \\\hline
			30 & $-8.83134288999026\cdot10^{-8}$ & $3.782652721710986\cdot10^{-7}$\\\hline
			31 & 0 & 0 \\\hline
			32 & $1.010610960272394\cdot10^{-7}$ & $2.418021923473667\cdot10^{-8}$\\\hline
			33 & 0 & 0 \\\hline
			34 & $6.606163280924149\cdot10^{-9}$ & $-2.689431432873997\cdot10^{-8}$\\\hline
			35 & 0 & 0 \\\hline
		\end{tabular}
	}
\end{table}

\begin{figure}[!]
	\centering\includegraphics[width=\textwidth,height=13cm]{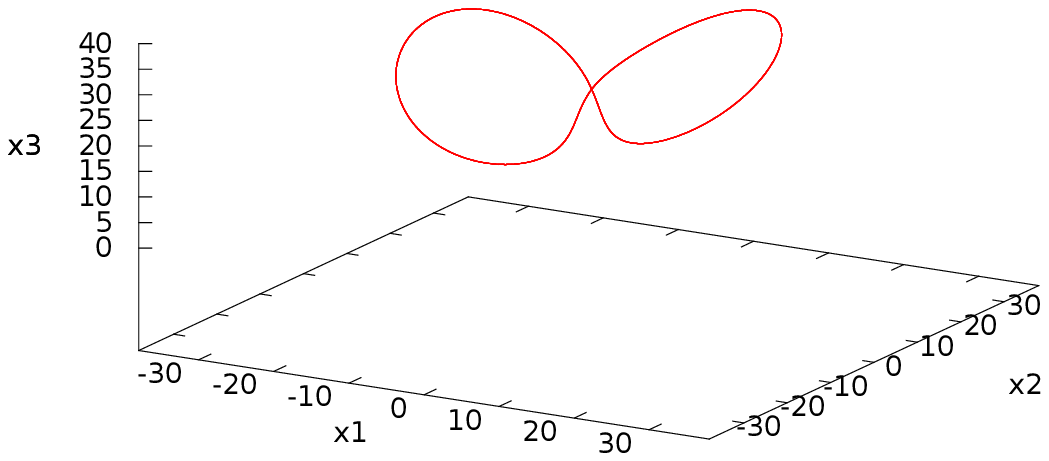}
	\caption{The cycle obtained by described method.}
\end{figure}

Tables 1--3 show the result of solving the system for $h=35$; the accuracy of the
Newton method is $10^{-8}$. The period value is obtained equal to $T=1.558652210$,
the initial condition for the obtained approximate periodic solution is
\begin{equation}\label{eq5}
\begin{array}{c}
\tilde{x}_1(0)=-2.147367631,\:\:\tilde{x}_2(0)=2.078048211,\:\:\tilde{x}_3(0)=27.
\end{array}
\end{equation}

The initial values \eqref{eq5} were checked on the period in a computer program that
implements the numerical integration of the system \eqref{eq1} by the modified power
series method \cite{8} with an accuracy of estimating the common term of the series
$10^{-25}$, 100 bits for mantissa real number and machine epsilon
$1.57772\cdot 10^{-30}$. With such parameters of the method, the approximate values
of the phase coordinates obtained by numerical integration were also verified by the
same numerical method, but in reverse time. The values in the reverse time coincide
with \eqref{eq5} up to the 9th character inclusive after the point. The resulting values
of $x_1(T)$, $x_2(T)$ and $x_3(T)$ coincide with \eqref {eq5} up to the 8th character
inclusive.

The cycle corresponding to \eqref{eq5} is shown in Fig. 1. Note that the cycle found
coincides with the first Viswanath cycle in \cite{6}, all signs after the point for
$T$ also coincide with the data from \cite{6}.

\section{Acknowledgements}

The reported study was funded by RFBR according to the research project 20-01-00347.

\end{document}